\documentclass[11pt,a4paper,final]{article}

\usepackage[latin1]{inputenc}

\usepackage[T1]{fontenc}
\usepackage[english]{babel}
\usepackage{amsmath,amssymb,amscd,latexsym}

\usepackage[draft]{showkeys}

\usepackage{theorem}
\theorembodyfont{\slshape}% penché plutot que italique

\numberwithin{equation}{section}
\newtheorem{Thm}[equation]{Theorem}

\newtheorem{Def}[equation]{Definition}
\newtheorem{Prop}[equation]{Proposition}
\newtheorem{Lem}[equation]{Lemma}
\newtheorem{Cor}[equation]{Corollary}

\theorembodyfont{\rmfamily}% en caracteres droit

\newtheorem{Rem}[equation]{Remark}

\newenvironment{proof}{\noindent{\sc Proof.}}{\hfill\qed}
\newenvironment{sproof}[1]{\noindent {\sc Proof of #1.}}{\hfill \qed}
\newcommand{\stab}{\operatorname{Stab}}

\newcommand{\N}{\mathbb{N}}
\newcommand{\Z}{\mathbb{Z}}

\newcommand{\qed}{\hfill $\square$}
\newcommand{\Pres}[2]{\left\langle{#1}\ \big\vert\ {#2}\right\rangle}

\title{Wreath products with the integers, proper actions and Hilbert space compression}
\author{Yves STALDER\footnote{Supported by the Swiss FNRS, grant  20-109130} \ and Alain VALETTE}
\date{March 29, 2006}

\begin{document}

\selectlanguage{english}

\maketitle

\begin{abstract}
We prove that the properties of acting metrically properly on some space with walls or some CAT(0) cube
complex are closed by taking the wreath product with $\Z$. We also give a lower bound for the (equivariant) Hilbert
space compression of $H\wr\Z$ in terms of the (equivariant) Hilbert space compression of $H$.
\end{abstract}

%---------------------------------------------------------------------------------------------------------------
\section*{Introduction}\label{Intro}

%All groups appearing in this note are supposed to be countable and discrete. See section \ref{Def} for
%definitions which do not appear in the introduction.

A {\it space with walls}, as defined by Haglund and Paulin \cite{HaPa}, is a pair $(X,\cal{W})$ where $X$ is
a set and $\cal{W}$ is a set of partitions of $X$ (called {\it walls}) into two classes, submitted to the
condition that any two points of $X$ are separated by finitely many walls.

The main examples of spaces with walls are given by $CAT(0)$ cube complexes (see \cite{BriH}), i.e. metric
polyhedral complexes in which each $k$-cell is isomorphic to the Euclidean cube $[-1/2,1/2]^k$, and the
gluing maps are isometries. Indeed, it is a result of Sageev \cite{Sag} that hyperplanes in a $CAT(0)$ cube
complex endow the set of vertices with a structure of space with walls (see \cite{ChatNib} and \cite{Nica}
for more on the relation between spaces with walls and $CAT(0)$ cube complexes).

Our first result is the following:

\begin{Thm}\label{PropWr}
Suppose that a group $H$ acts metrically properly either on some on some space with walls, or on some
$CAT(0)$ cube complex. Then the wreath product $H\wr\Z\,:=\,(\bigoplus_{\Z} H)\rtimes\Z$ satisfies the same
property.
\end{Thm}

Guentner-Kaminker defined the \emph{Hilbert space compression} and the \emph{equivariant Hilbert space
compression} for any unbounded metric space (endowed with a group action in the second case) \cite{GuKa}.
Since we will deal with uniformly discrete\footnote{That
is, there exists a constant $\delta > 0$ such that $d(x,y) \geqslant \delta$ whenever $x\neq y$.} spaces, the
following definitions are equivalent to theirs.

Let $(X,d)$ be a uniformly discrete metric space. We define the \emph{Hilbert space
compression} of $X$ as the supremum of the numbers $\alpha \in [0,1]$ such that there exists a Hilbert space
$\mathcal{H}$, positive constants $C_1,C_2$ and a map $f:X\to \mathcal{H}$ with
\[
C_1 \cdot d(x,y)^\alpha \leqslant || f(x) - f(y) || \leqslant C_2 \cdot d(x,y) \ \forall x,y \in X \ .
\]
It is denoted by $R(X,d)$ and it is a quasi-isometry invariant of $(X,d)$.
If $H$ is a group acting on $(X,d)$ by isometries, the \emph{equivariant Hilbert space
compression} of $X$ is the supremum of the numbers $\alpha \in [0,1]$ such that there exists a Hilbert
space $\mathcal{H}$ endowed with an action of $H$ by affine isometries, positive constants $C_1,C_2$ and
a $H$-equivariant map $f:X\to \mathcal{H}$ with
\[
C_1 \cdot d(x,y)^\alpha \leqslant || f(x) - f(y) || \leqslant C_2 \cdot d(x,y) \ \forall x,y \in X \ .
\]
It is denoted by $R_H(X,d)$. One has trivially $R_H(X,d) \leqslant R(X,d)$.

We may view a group $H$ as a metric space thanks to the word length associated with some (not necessarily
finite) generating subset $S$. We denote then by $R(H,S)$ the Hilbert space compression and by $R_H(H,S)$ the
equivariant Hilbert space compression. In case $H$ is finitely generated, note that,
up to bilipschitz equivalence, the word metric does not depend on the finite generating set,
so that the compressions do not depend on the choice of a finite generating set. In this case
we write $R(H)$ and $R_H(H)$ for the corresponding compressions. We also use these shorter notations in the
general case if there is no ambiguity about the generating set. It is a remarkable observation of Gromov (see
\cite[Proposition 4.4]{CTV} for a proof) that $R(H) = R_H(H)$ for $H$ finitely generated and amenable.

The first examples of finitely generated groups whose Hilbert space compression is different from $0$ and $1$
appeared recently in \cite{AGS}: Thompson's group
$F$ and the wreath product $\Z\wr\Z$ (see the end of Section \ref{Hcomp} for more on this).
Our next Theorem allows in particular to construct more examples.

Given a generating set $S$ for $H$, if $\Gamma=H\wr\Z$, we always take
$\Sigma= S\cup\{s\}$ as generating set for $\Gamma$, where $s$ is the positive generator of $\Z$.

\begin{Thm}\label{HilbC}
Let H be a group, with generating set $S$ and let $\Gamma= H\wr \Z$. The non equivariant and equivariant
Hilbert space compressions satisfy:
\begin{align*}
  R(H,S) \geqslant R(\Gamma,\Sigma) & \geqslant \frac{R(H,S)}{R(H,S) + 1} \ ;  \\
  R_H(H,S) \geqslant R_\Gamma(\Gamma,\Sigma) & \geqslant \max\left\{ R_H(H,S) - \frac 12 \, , \,
  \frac{R_H(H,S)}{2R_H(H,S) + 1} \right\} \ .
\end{align*}
\end{Thm}
In order to select the best bound, we mention that one has $t - 1/2 \geqslant t/(2t+1)$ if and only if
$t\geqslant (1+ \sqrt 5)/4 \cong 0.809...$ (for $t\in [0,1]$). Gromov's remark gives immediately a stronger
estimate for the equivariant compression.
\begin{Cor}\label{HilbCmoy}
Let H be a finitely generated and amenable group and let $\Gamma= H\wr \Z$. The equivariant Hilbert space
compression satisfies:
\[
R_H(H) \geqslant R_\Gamma(\Gamma) \geqslant \frac{R_H(H)}{R_H(H) + 1} \ .
\]
\end{Cor}

The proofs of Theorems \ref{PropWr} and \ref{HilbC} rest on a similar idea: we express $H\wr\Z$ as an HNN-extension in two
different ways, which provide two different actions of $H\wr\Z$ on a tree. In Theorem \ref{PropWr} we use the product of these two
trees, while in Theorem \ref{HilbC} we appeal to the affine actions naturally associated with each of these trees (see section
7.4.1 in \cite{CCJJV}).

%---------------------------------------------------------------------------------------------------------------
\section{Preliminaries: wreath products and trees}\label{Def}

%\paragraph{HNN-extensions and Bass-Serre trees.}
Let $\Lambda$ be a group, $H$ a subgroup and
$\vartheta: H \to \Lambda$ an injective homomorphism. The \emph{HNN-extension} with \emph{basis} $\Lambda$ and
\emph{stable letter} $t$ relatively to $H$ and $\vartheta$ is defined by $HNN(\Lambda,H,\vartheta) =
\Pres{\Lambda, t}{t^{-1}ht = \vartheta(h) \ \forall h \in H}$.

Our definition of graphs and trees are those of \cite{Ser}.
Given an HNN-extension $\Gamma = HNN(\Lambda,H,\vartheta)$, the associated \emph{Bass-Serre tree} is defined by
\[
V(T) = \Gamma/\Lambda \ ; \ E(T) = \Gamma/H \sqcup \Gamma/\vartheta(H) \ ; \ \overline{\gamma H} = \gamma t \vartheta(H) \
; \ \overline{\gamma \vartheta(H)} = \gamma t^{-1} H \ ;
\]
\[
(\gamma H)^- = \gamma \Lambda \ ; \ (\gamma H)^+ = \gamma t \Lambda \ ; \
(\gamma \vartheta(H))^- = \gamma \Lambda \ ; \ (\gamma \vartheta(H))^+ = \gamma t^{-1} \Lambda
\]
where, given an edge $e$, its origin is denoted by $e^-$ and its terminal vertex by $e^+$. It is a tree
\cite[Theorem 12]{Ser}. We turn $T$ to an oriented tree by setting $Ar_+(T) = \Gamma/H$, $(\gamma H)^- =
\gamma \Lambda$, $(\gamma H)^+ = \gamma t \Lambda$ and the $\Gamma$-action on $T$ preserves this
orientation. Moreover we remark that the oriented tree is bi-regular: for each vertex of $T$ the outgoing
edges are in bijection with $\Lambda / H$ and the incoming edges are in bijection with $\Lambda /
\vartheta(H)$.

We turn to wreath products. Let $G,H$ be groups. We set
\[
\Lambda = H^{(G)} = \bigoplus_{g\in G} H = \{ \lambda:G \to H \text{ with finite support } \} \ .
\]
The group $G$ acts on $\Lambda$ by automorphisms: $(g\cdot\lambda)(x) = \lambda(g^{-1} x)$. The
\emph{wreath product} $H\wr G$ is the semi-direct product $\Lambda \rtimes G$, with respect to the action above.
The group $H$ embeds in $H\wr G$ as the copy of index $1_G$. It is easy to see that, given generating sets of
$G$ and $H$, their union generates $H\wr G$.

In case $G=\Z$, one may express $H\wr\Z$ as an HNN-extension in two ways (we denote by $s$ the generator 
of $\Z$ in $H\wr\Z$ and by $t_+,t_-$ the stable letters of the HNN-extensions)\footnote{
  The most common way is probably the following: set $\vartheta : \Lambda \to \Lambda ; 
  \vartheta(\lambda)_n = \lambda_{n-1}$. One has
  $HNN(\Lambda,\Lambda,\vartheta) = H\wr\Z$ and the isomorphism is given by $\lambda \mapsto \lambda$
  and $t\mapsto s^{-1}$. Nevertheless, this expression will be useless in this article.}:
\begin{enumerate}
  \item Set $\Lambda_+ = \bigoplus_{n\geqslant 0}H$ and $\vartheta_+: \Lambda_+ \to \Lambda_+$ given by
  $\vartheta_+(\lambda)_0=1_H$ and $\vartheta_+(\lambda)_n= \lambda_{n-1}$ for $n\geqslant 1$. One has
  $HNN(\Lambda_+,\Lambda_+,\vartheta_+)=H\wr\Z$ and the isomorphism is given by $\lambda \mapsto \lambda$
  and $t_+\mapsto s^{-1}$;

  \item Set $\Lambda_- = \bigoplus_{n\leqslant 0}H$ and $\vartheta_-: \Lambda_- \to \Lambda_-$ given by
  $\vartheta_-(\lambda)_0=1_H$ and $\vartheta_-(\lambda)_n= \lambda_{n+1}$ for $n\leqslant -1$. One has
  $HNN(\Lambda_-,\Lambda_-,\vartheta_-)=H\wr\Z$ and the isomorphism is given by $\lambda \mapsto \lambda$
  and $t_-\mapsto s$;
\end{enumerate}

Given a wreath product $H\wr\Z$, we will denote by $T_+$, respectively $T_-$, the Bass-Serre tree associated
to the second, respectively third, HNN-extension above. We take as base points (when necessary) the vertices
$\Lambda_+$ and $\Lambda_-$.

%\paragraph{CAT(0) cube complexes and spaces with measured walls.} A \emph{cube complex} is a metric
%polyhedral complex in which each $k$-cell is isomorphic to the euclidean cube $[-1/2,1/2]^k$, and
%the gluing maps are isometries.

%We refer to \cite{BriH} for the definition of the CAT(0) condition and to \cite{HaPa} for the definition
%of a space with walls.

%\paragraph{Hilbert space compressions}

%--------------------------------------------------------------------------------------------------
%\section{Observations}\label{Obs}

We collect now some observations about the $H\wr\Z$-actions on $T_+$ and $T_-$ which will be relevant
in the next sections. Set $\Gamma = H\wr\Z$ and $\gamma = (\lambda,n) \in \Gamma$. If $\lambda$ is nontrivial, we
set $m= \min\{ k\in\Z : \lambda_k \neq 1_H \}$ and  $M= \max\{ k\in\Z : \lambda_k \neq 1_H \}$.

\begin{Lem}\label{obs1}
If $\lambda = 1$, one has $d_{T_+}(\Lambda_+, \gamma \Lambda_+) = |n| = d_{T_-}(\Lambda_-, \gamma \Lambda_-)$.
\end{Lem}
The proof is obvious.

\begin{Lem}\label{obs2}
If $\lambda \neq 1$, the distances $d_{T_\pm}(\Lambda_\pm, \gamma \Lambda_\pm)$ are given by formulas:
  \begin{eqnarray*}
    d_{T_+}(\Lambda_+, \gamma \Lambda_+) & = & \left\{
    \begin{array}{cl}
      |n| & \text{ if } n \leqslant m \text{ or } m \geqslant 0   \\
      n-2m & \text{ if } n > m \text{ and } m < 0  \\
    \end{array}
    \right. \ ; \\
    d_{T_-}(\Lambda_-, \gamma \Lambda_-) & = & \left\{
   \begin{array}{cl}
      |n| & \text{ if } n \geqslant M \text{ or } M \leqslant 0   \\
      2M-n & \text{ if } n < M \text{ and } M > 0  \\
    \end{array}
    \right. \ . \\
  \end{eqnarray*}
In particular, the inequalities $d_{T_+}(\Lambda_+, \gamma \Lambda_+) \geqslant -m$, $d_{T_-}(\Lambda_-, \gamma \Lambda_-)
\geqslant M$, $d_{T_+}(\Lambda_+, \gamma \Lambda_+) \geqslant |n|$ and $d_{T_-}(\Lambda_-, \gamma \Lambda_-) \geqslant |n|$
hold.
\end{Lem}

\begin{proof}
We prove the first equality, leaving the second one, which is very similar, to the reader. We remark
that $\gamma = \lambda s^n$ and that, for any $k\in Z$, the stabilizer of the vertex $t_+^k\Lambda_+$
satisfies
\begin{equation}\label{eqstab}
\stab(t_+^k\Lambda_+) = t_+^k\Lambda_+ t_+^{-k} = s^{-k}\Lambda_+ s^k = \bigoplus_{i\geqslant -k}H \ .
\end{equation}
Suppose first that $m\geqslant 0$. Then $\lambda$ stabilizes the vertex $\Lambda_+$, so that we get
$d(\Lambda_+, \gamma\Lambda_+) = d(\lambda \Lambda_+, \lambda s^n\Lambda_+) = d(\Lambda_+, s^n \Lambda_+) = |n|$.
If $n\leqslant m$, the vertex $s^n\Lambda_+ = t_+^{-n}\Lambda_+$ is stabilized by $\lambda$, so that
$d(\Lambda_+, \gamma\Lambda_+) = d(\Lambda_+, s^n \Lambda_+) = |n|$.

It remains to treat the case $n > m$ and $m < 0$. The vertices on the geodesic from $t_+^{-m} \Lambda_+$ to
$t_+^{-n} \Lambda_+$ are $t_+^{-m}\Lambda_+, t_+^{-m-1}\Lambda_+, \ldots, t_+^{-n}\Lambda_+$. By (\ref{eqstab}), the
vertex $t_+^{-m} \Lambda_+$ is stabilized by $\lambda$ and $t_+^{-m-1} \Lambda_+$ is not. Thus, the geodesic from
$\Lambda_+$ to $\gamma \Lambda_+$ passes through $t_+^{-m} \Lambda_+$, so that we get
\[
d(\Lambda_+, \gamma\Lambda_+) = d(\Lambda_+, t_+^{-m}\Lambda_+) + d(t_+^{-m}\Lambda_+, \gamma\Lambda_+) = -m + n-m = n-2m \ .
\]
The proof is complete.
\end{proof}

Let us now state a formula computing the length of an element of $H\wr\Z$, which
is a direct consequence of \cite[Theorem 1.2]{Par}. Note that, even if the theorem was stated for finitely
generated groups, it also applies in our case.

\begin{Prop}\label{ParticPar} Keep the above notations.
Let $\gamma = (\lambda,n) \in \Gamma = H\wr\Z$. In case $\lambda = 1$, one has
$|\gamma| = |n|$, while in case $\lambda \neq 1$, the length of $\gamma$ satisfies:
\[
|\gamma| =  L_\Z(\gamma) + \sum_{i\in\Z} |\lambda_i| \ ,
\]
where $L_\Z(\gamma)$ denotes the length of the shortest path starting from $0$, ending at $n$ and passing
through $m$ and $M$ in the (canonical) Cayley graph of $\Z$.
\end{Prop}

The length $L_\Z(\gamma)$ appearing in Proposition \ref{ParticPar} can be estimated as follows:

\begin{Prop}\label{EstimL}
Let $\gamma\in \Gamma = H \wr \Z$. The following inequalities hold:
\[
d_{T_\pm}(\Lambda_\pm,\gamma\Lambda_\pm) \leqslant L_\Z(\gamma) \leqslant
d_{T_+}(\Lambda_+,\gamma\Lambda_+) + d_{T_-}(\Lambda_-,\gamma\Lambda_-)
\]
\end{Prop}

\begin{proof}
If $\gamma=(1,n)$, the result is obvious. We suppose now $\gamma=(\lambda,n)$ with $\lambda \neq 1$.
The proof is then a distinction of eight cases which are listed in the following tabular:
\[
\begin{array}{c|c|c||c|c|c|c}
  % after \\ : \hline or \cline{col1-col2} \cline{col3-col4} ...
    n         &     m       &       M     & d_{T_+} & d_{T_-}   & L_\Z      & d_{T_+} + d_{T_-}     \\  \hline
  \geqslant 0 & \geqslant 0 &       >n    &   n     & 2M-n      & 2M-n      & 2M                    \\
  \geqslant 0 & \geqslant 0 & \leqslant n &   n     &   n       &   n       & 2n                    \\
  \geqslant 0 &      < 0    &       >n    &  n-2m   & 2M-n      & 2M-2m-n   & 2M - 2m               \\
  \geqslant 0 &      < 0    & \leqslant n &  n-2m   &    n      &   n-2m    & 2n - 2m               \\ \hline
       < 0    &      < n    & \leqslant 0 &  n-2m   &   -n      & n -2m     & -2m                    \\
       < 0    & \geqslant n & \leqslant 0 &   -n    &    -n     &   -n      & -2n                    \\
       < 0    &      < n    &     > 0     &  n-2m   &   2M-n    & 2M-2m+n   & 2M-2m                  \\
       < 0    & \geqslant n &     > 0     &   -n    &   2M-n    &   2M-n    & 2M-2n                  \\ \hline
\end{array}
\]
The values of $d_{T_\pm}(\Lambda_\pm,\gamma\Lambda_\pm)$ come from Lemma \ref{obs2}; those of $L_\Z(\gamma)$
are easy to compute. We now observe that the result is true in the eight cases.
\end{proof}

Combining Propositions \ref{ParticPar} and \ref{EstimL}, one obtains immediately:

\begin{Cor}\label{Estim}
Let $\gamma = (\lambda,n) \in \Gamma = H\wr\Z$. The following inequalities hold:
\[
d_{T_\pm}(\Lambda_\pm,\gamma\Lambda_\pm) + \sum_{i\in\Z} |\lambda_i| \leqslant|\gamma| \leqslant
d_{T_+}(\Lambda_+,\gamma\Lambda_+) + d_{T_-}(\Lambda_-,\gamma\Lambda_-) + \sum_{i\in\Z} |\lambda_i|
\]
\end{Cor}

%---------------------------------------------------------------------------------------------------
\section{Metrically proper actions}\label{PropAct}

Let us consider a group $G$, acting by isometries on a metric space $X$.
\begin{Def}
The action is \emph{metrically proper} if, whenever $B$ is a bounded subset of $X$, the set
$\{ g\in G : g\cdot B \cap B \neq \varnothing \}$ is finite.
\end{Def}
From now on, we shall write ``proper'' instead of ``metrically proper''. Let us recall that
the action is proper if and only if the following property holds, for some $z\in X$:
\begin{equation*}\label{locProp}
\text{for any } R>0, \text{ the set } \{ g\in G : d(z,g\cdot z)\leqslant R \} \text{ is finite.} \tag{$Prop_z$}
\end{equation*}

Let now $\mathcal{F} = (X_i,b_i)_{i\in I}$ be a family of pointed metric spaces and let $p \geqslant 1$.
We call \emph{$\ell^p$-product}  of the family the space
\[
\ell^p(\mathcal{F}) = \left\{ x\in \prod_{i\in I} X_i : \sum_{i\in I} d_i(b_i,x_i)^p
< +\infty \right\} \ .
\]
It is a metric space with metric $\delta(x,y) = \left(\sum_{i\in I} d_i(x_i,y_i)^p
\right)^{1/p}$. We set $(b_i)_{i\in I}$ as base point. Consider now the case $(X_i,b_i)=(X,b)$ for all
$i\in \Z$. One has:
\[
\ell^p(I;X,b):= \ell^p(\mathcal{F}) = \left\{ \phi:I\to X : \sum_{i\in I} d(b,\phi_i)^p < +\infty \right\} \ .
\]
If a group $H$ acts by isometries on $X$, the group $H\wr G$ acts by isometries on $\ell^p(G;X,b)$
in the following way:
\begin{equation}\label{ActWr}
\left\{
  \begin{array}{cccl}
  (\lambda \cdot \phi)_g & = & \lambda_g \cdot \phi_g & \text{ for } \lambda\in \bigoplus_{g\in G}H \ ; \\
  (g\cdot\phi)_{g'}         & = & \phi_{g^{-1}g'}             & \text{ for } g \in G \ .
  \end{array}
\right.
\end{equation}

Given $G$ infinite, observe that, even if the action of $H$ is proper, the action of $H \wr G$ on $\ell^p(G;X,b)$
is \textit{not}. Indeed, there is a $G$-globally fixed point on $\ell^p(G;X,b)$.

Theorem \ref{PropWr} will follow from the following statement.

\begin{Prop}\label{PropGen}
Let $H$ be a group acting properly on a metric space $X$, $b\in X$ and $p\geqslant 1$. Then, the action of
$\Gamma= H\wr\Z$ on $T_+ \times T_- \times \ell^p(\Z;X,b)$, where the product is endowed with the $\ell^p$
metric, is proper.
\end{Prop}

\begin{proof}
We are going to prove property (\ref{locProp}) for $z=(\Lambda_+,\Lambda_-,(b)_{i\in\Z})$. Thus let
$R>0$ and $A= \{ \gamma \in \Gamma : d(z,\gamma\cdot z)\leqslant R \}$. Take $\gamma=(\lambda,n) \in A$.
We have $d_{T_+}(\Lambda_+, \gamma \Lambda_+) \leqslant R$, $d_{T_-}(\Lambda_-, \gamma \Lambda_-) \leqslant R$ and
$\sum_{i\in\Z} d(b,\lambda_i \cdot b)^p \leqslant R^p$.

By lemmata \ref{obs1} and \ref{obs2}, one has $M\leqslant R$, $m\geqslant -R$ (if $M$ and $m$ are defined)
and $|n|\leqslant R$. Set $B= \{ h \in H : d(b,h\cdot b) \leqslant R \}$. It is a finite set since the
$H$-action is proper.

Hence, one has $|n| \leqslant R$, $\lambda_i = 1_H$ for $|i| > R$ and $\lambda_i \in B$ for $|i|\leqslant R$.
This leaves finitely many choices for $\gamma$, and proves thus that $A$ is finite.
\end{proof}

\begin{Rem}
The space $T_+ \times T_- \times \ell^p(\Z;X,b)$ is canonically isometric to the product $\ell^p(\mathcal{F})$
with $I = \{+,-\} \cup \Z$ and $\mathcal{F}$ given by $\mathcal{F}(+) = (T_+,\Lambda_+)$,
$\mathcal{F}(-) = (T_-,\Lambda_-)$ and $\mathcal{F}(i) = (X,b)$ for $i\in\Z$.
\end{Rem}

\begin{sproof}{of Theorem \ref{PropWr}}
We recall first that a tree is a CAT(0) cube complex, hence a space with walls.

It is shown in \cite[Section 5]{CMV} that a $\ell^1$-product of spaces with measured walls carries the same structure.
Moreover, we remark that, particularizing the construction to spaces with walls, one gets a space with walls.
Hence, we get the conclusion for spaces with walls by proposition \ref{PropGen}.

Given a CAT(0) cube complex $Y$, we denote by $Y^{(k)}$ the set of $k$-cells in $Y$.
Take now a family $\mathcal{F}=(X_i,b_i)_{i\in I}$ of CAT(0) cube complexes with $b_i \in X_i^{(0)}$
and set $\mathcal{F}^{(0)} = (X^{(0)}_i,b_i)_{i\in I}$ for $k\in \N$. We are going to construct a subspace
$X$ of $\ell^2(\mathcal{F})$ which is a CAT(0) cube complex.

We define first $X^{(0)} = \ell^2(\mathcal{F}^{(0)})$. Since the distance between two distinct vertices is
at least $1$, one has
\[
X^{(0)} = \bigoplus_{i\in I} (X_i^{(0)},b_i) := \left\{ v \in \prod_{i\in I} X_i^{(0)} :
\{ i\in I : v_i \neq b_i \} \text{ is finite } \right\} \ .
\]
For $k \geqslant 1$, we define then the set of $k$-cells as
\[
X^{(k)} = \left\{
\begin{array}{c}
  c \in \prod_{i\in I} (X_i^{(0)} \cup \ldots \cup X_i^{(k)} ) : \\
  \sum_{i\in I} \dim(c_i) = k \text{ and } \{ i\in I : c_i \neq b_i \} \text{ is finite }
\end{array}
\right\} \ .
\]
It is clear that every $k$-cell, as a subset of $\ell^2(\mathcal{F})$, is isometric to $[-1/2,1/2]^k$. If
$c\in X^{(k)}$, the faces of $c$ are the $(k-1)$-cells $c'$ such that $c'_j$ is a face of $c_j$ for some $j$
and $c'_i=c_i$ for $i \neq j$. The gluing maps are isometric. Finally, the space $\ell^2(\mathcal{F})$
inherits the CAT(0) property, so that $X$ is a CAT(0) cube complex.

Suppose now that H acts on a CAT(0) cube complex $Y$ and take $v_0$ a vertex of $Y$.
We consider the family $\mathcal{F}$ given by $I= \{+,-\} \cup \Z$, $\mathcal{F}(+) = (T_+,\Lambda_+)$,
$\mathcal{F}(-) = (T_-,\Lambda_-)$ and $\mathcal{F}(i) = (Y,v_0)$ for $i\in\Z$. The action of $H\wr \Z$
on $\ell^2(\mathcal{F})$ is proper by proposition \ref{PropGen} and the CAT(0) cube complex $X$ constructed
as above is an invariant subset, so that it is endowed with a proper action of $H\wr \Z$ too.
\end{sproof}

\begin{Rem}\label{StabHaa1}
The same techniques can be used, if $H$ acts properly on some Hilbert space  $\mathcal{H}$, to prove that
$H\wr\Z$ acts properly on the Hilbert direct sum
$\ell^2(E(T_+)) \oplus \ell^2(E(T_-)) \oplus \bigoplus_{i\in\Z} \mathcal{H}$. Hence,
we recover the known fact that Haagerup property is preserved by taking wreath products with $\Z$
\cite[Proposition 6.1.1 and Example 6.1.6]{CCJJV}. The interest of our technique is that we obtain an
\emph{explicit} proper action of $H\wr\Z$, knowing a proper action of $H$.
\end{Rem}

\begin{Rem}\label{StabHaa2}
It is known \cite[Theorem 1]{CMV}, that a discrete group satisfies the Haagerup property if and only if it
acts properly on some space with measured walls. It follows from Remark \ref{StabHaa1} that whenever $H$
acts properly on a space with measured walls, the same holds for $H\wr\Z$. Again, our techniques give
an explicit action, as Theorem \ref{PropWr} is also valid for spaces with measured walls.
\end{Rem}

%---------------------------------------------------------------------------------------------------
\section{Hilbert space compression: Theorem \ref{HilbC}}\label{Hcomp}

We recall that a map $f: X \to Y$
between metric spaces is \emph{Lipschitz} if there exists $C>0$ such that $d_Y(f(x),f(y)) \leqslant C \cdot
d_X(x,y)$ for all $x,y \in X$. Given a Lipschitz map $f:X \to \mathcal H$ (whose range is a Hilbert space),
we set\footnote{It does not coincide with the
\emph{asymptotic compression} of $f$ defined in \cite{GuKa}.} 
$R_f$ to be the supremum of the numbers $\alpha \in [0,1]$ such that there exists $D>0$ with $D \cdot
d_X(x,y)^\alpha \leqslant ||f(x) - f(y)||$ for all $x,y\in X$.

Given a generating set $S$ of a group $H$, we
recall our convention to take $\Sigma= S\cup\{s\}$ as generating set for $H\wr\Z$, where $s$ is the positive
generator of $\Z$. In order to simplify notations, we do not mention explicitly $S$ and $\Sigma$ in this section.

The goal of this section is to prove Theorem \ref{HilbC}. The key result in this way is the following:

\begin{Prop}\label{HilbWr}
Let $H$ be a group (with a generating set $S$) and $\Gamma = H\wr\Z$. Suppose that maps ${f: H \to \mathcal{H}}$ and
$f_\pm : V(T_\pm) \to \mathcal{H}_\pm$ are Lipschitz with ${R_{f_+} = R_{f_-} >0}$ and
${R_f>0}$. Then consider the map
\[
\sigma: \Gamma \to \mathcal{H}' :=  \mathcal{H}_+ \oplus \mathcal{H}_- \oplus \bigoplus_{i\in\Z} \mathcal{H} \ ,
\]
where, given $\gamma=(\lambda,n) \in H\wr\Z$, we set $\sigma(\gamma)_\pm = f_\pm(\gamma \Lambda_\pm)$ and
$\sigma(\gamma)_i = f(\lambda_i)$ for $i\in \Z$. It satisfies $R_\sigma \geqslant R_f\cdot R_{f_\pm}/(R_f+R_{f_\pm})$
and $R_\sigma \geqslant \min\{R_{f_\pm}, R_f - \frac 12\}$

Moreover, if $f$ is $H$-equivariant and if $f_\pm$ are $\Gamma$-equivariant (with respect to
some actions by affine isometries), then there exists a $H$-action by affine isometries on $\mathcal{H}'$
such that $\sigma$ is $\Gamma$-equivariant.
\end{Prop}

\begin{proof}
We show first that $\sigma$ is Lipschitz (the reader could remark that it is trivial if $H$ is
a finitely generated group; however, this case is also covered by the proof below). Let us take
$C,C_+,C_- >0$  such that
\begin{eqnarray*}
  ||f(s) - f(t)|| & \leqslant & C \cdot |s^{-1}t|  \quad \forall s,t \in H \ ; \\
  ||f_\pm(u) - f_\pm (v)|| & \leqslant & C_\pm \cdot d_{T_\pm}(u,v) \quad \forall u,v \in V(T_\pm) \ .
\end{eqnarray*}
Let $x,y \in \Gamma$. We set $\gamma = x^{-1} y$ and  write $x = (\xi,p)$, $y=(\eta,q)$, $\gamma=(\lambda,n)$ in
$H\wr\Z = \Lambda \rtimes \Z$ (so that $n=q-p$ and $\lambda_i = \xi_{i-p}^{-1} \eta_{i-p}$).
One has then
\[
 \left(\sum_{i\in\Z} || \sigma(x)_i - \sigma(y)_i ||^2 \right)^{\frac 12} \leqslant \sum_{i\in\Z} || f(\xi_i) - f(\eta_i) ||
    \leqslant \sum_{j\in\Z} C\cdot |\lambda_j| \leqslant C\cdot |\gamma| \ .
\]
Moreover, using Corollary \ref{Estim} for the last step, it comes:
\[
|| \sigma(x)_\pm - \sigma(y)_\pm || \leqslant C_\pm \cdot d_{T_\pm}(x\Lambda_\pm, y\Lambda_\pm)
        = C_\pm \cdot d_{T_\pm}(\Lambda_\pm, \gamma\Lambda_\pm) \leqslant C_\pm \cdot |\gamma| \ .
\]
Thus, we get finally $|| \sigma(x) - \sigma(y) || \leqslant (C_+ + C_- + C) \cdot |x^{-1}y|$,
which proves that $\sigma$ is Lipschitz, as desired.

We now turn to the estimation of $R_\sigma$, Fix any $\alpha,\beta$ such that $0<\alpha<R_f$ and
$0<\beta<R_{f_\pm}$. There exists constants $C,C_+,C_- >0$ such that:
\begin{eqnarray*}
  ||f(s) - f(t)|| & \geqslant & C \cdot |s^{-1}t|^\alpha \quad \forall s,t \in H \ ; \\
  ||f_\pm(u) - f_\pm (v)|| & \geqslant & C_\pm \cdot d_{T_\pm}(u,v)^\beta \quad \forall u,v \in V(T_\pm) \ .
\end{eqnarray*}

We notice first that $\sigma$ is injective. More precisely, for any $x,y\in\Gamma$, one has
\begin{equation}\label{UnifInj}
x\neq y \ \Longrightarrow \ ||\sigma(x)- \sigma(y)|| \geqslant \min\{C,C_+,C_-\} \ .
\end{equation}
Indeed, we express $x=(\xi,p)$ and $y=(\eta,q)$ as above. If $p\neq q$, we obtain $||\sigma(x)_\pm- \sigma(y)_\pm|| \geqslant
C_\pm \cdot d_{T_\pm}(x\Lambda_\pm,y\Lambda_\pm)^\beta \geqslant C_\pm$ and if $\xi_i \neq \eta_i$ for some
$i$, we obtain $||\sigma(x)_{i} - \sigma(y)_{i}|| \geqslant C\cdot |\xi_i^{-1}\eta_i|^\alpha \geqslant C$.

% Let us take any $\zeta \in \ ]0,1[$ and some $a\in \ ]0,1[$ such
%that $a^{1-\zeta} > 1/2$.

Let us take $x$, $y$ and $\gamma$ as above.
According to Corollary \ref{Estim}, one (at least) of the following cases occurs. We
treat them separately.
As the case $x=y$ is trivial, we assume $x\neq y$, that is $|\gamma| \geqslant 1$, in what follows.
\begin{description}
  \item[(a) Case $d_{T_+}(\Lambda_+,\gamma\Lambda_+) \geqslant \frac 13 |\gamma|$:] We obtain
    \begin{eqnarray*}
    || \sigma(x) - \sigma(y) || & \geqslant & || \sigma(x)_+ - \sigma(y)_+ || \geqslant
            C_+ \cdot d_{T_+}(x\Lambda_+,y\Lambda_+)^{\beta} \\
     & = & C_+ \cdot d_{T_+}(\Lambda_+,\gamma\Lambda_+)^{\beta} \geqslant
     \frac{C_+}{3^\beta} |x^{-1}y|^{\beta}  \ .
    \end{eqnarray*}
  \item[(b) Case $d_{T_-}(\Lambda_-,\gamma\Lambda_-) \geqslant \frac 13 |\gamma|$:] We obtain the same way
    \[
    || \sigma(x) - \sigma(y) ||  \geqslant  || \sigma(x)_- - \sigma(y)_- || \geqslant
    \frac{C_-}{3^\beta} |x^{-1}y|^{\beta}  \ .
    \]
  \item[(c) Case $\sum_{i\in\Z} |\lambda_i| \geqslant \frac 13 |\gamma|$:] We establish two
    independant estimates.

    First, for all $i \in \Z$, one has $|| \sigma(x)_i - \sigma(y)_i || = || f(\xi_i) - f(\eta_i) ||
    \geqslant C |\xi_i^{-1} \eta_i|^\alpha = C |\lambda_{i+p}|^\alpha$. Hence, using Cauchy-Schwarz inequality
    for the third step below and $\alpha \leqslant 1$ for the fourth one, we obtain
    \begin{align*}
    || \sigma(x) - \sigma(y) || & \geqslant \left( \sum_{i\in \Z} || \sigma(x)_i - \sigma(y)_i ||^2 \right) ^{\frac 12}
                                        \geqslant C \left( \sum_{j\in \Z} |\lambda_j|^{2\alpha} \right) ^{\frac 12}       \\
                                & \geqslant \frac{C}{\sqrt{M-m+1}} \sum_{j=m}^M |\lambda_j|^{\alpha}
                                        \geqslant \frac{C}{\sqrt{M-m+1}} \left( \sum_{j=m}^M |\lambda_j| \right)^{\alpha}
    \end{align*}
    By Proposition \ref{ParticPar}, one has $|\gamma| \geqslant M-m+1$, so that we obtain
    \begin{equation*}\label{eq1Hwr} \tag{$\ast$}
    || \sigma(x) - \sigma(y) || \geqslant  \frac{C}{\sqrt{|\gamma|}} \left( \frac 13 |\gamma| \right)^{\alpha}
    \geqslant \frac{C}{3^\alpha} |\gamma|^{\alpha - \frac 12} \ .
    \end{equation*}
    This is our first estimate for case (c).

    Second, we fix any $\zeta \in \ ]0,1[ \ $. Then, either there exists $k \in \Z$ such that
    ${|\lambda_k| \geqslant(\frac 13 \cdot |\gamma|)^\zeta}$, or
    one has $M-m+1 \geqslant (\frac 13 \cdot |\gamma|)^{1-\zeta}$. We distinguish the two subcases:
    \begin{itemize}
        \item if there exists $k \in \Z$ such that $|\lambda_k| \geqslant (\frac 13 \cdot |\gamma|)^\zeta$, we have
        \begin{align*}
        || \sigma(x) - \sigma(y) || & \geqslant || \sigma(x)_{k-p} - \sigma(y)_{k-p} ||
                                        \geqslant C \cdot |\xi_{k-p}^{-1} \eta_{k-p}|^\alpha \\
                                    & = C \cdot |\lambda_k|^{\alpha} \geqslant \frac C{3^{\alpha\zeta}}
                                        |\gamma|^{\alpha\zeta} \ ;
        \end{align*}

        \item in case $M-m+1 \geqslant (\frac 13 \cdot |\gamma|)^{1-\zeta}$, having $L_\Z(\gamma)
        \geqslant M-m$ by definition, Proposition \ref{EstimL} gives
        \[
        d_{T_+}(\Lambda_+,\gamma\Lambda_+) + d_{T_-}(\Lambda_-,\gamma\Lambda_-) \geqslant L_\Z(\gamma)
        \geqslant \left( \frac 13 \cdot |\gamma| \right)^{1-\zeta} - 1 \ .
        \]
        Thus, $\exists \ s\in\{+,-\}$ such that $d_{T_s}(\Lambda_s,\gamma\Lambda_s) \geqslant
        \frac 12 (\frac 13 \cdot |\gamma|)^{1-\zeta} - \frac 12$. For $|\gamma| \geqslant 4$, there exists
        $K>0$ such that $d_{T_s}(\Lambda_s,\gamma\Lambda_s) \geqslant K \cdot |\gamma|^{1-\zeta}$,
        so that $||\sigma(x) - \sigma(y)|| \geqslant C_s K^\beta \cdot |x^{-1}y|^{\beta(1-\zeta)}$
        as in cases (a)-(b).

        Otherwise, for $|\gamma| \leqslant 3$, equation (\ref{UnifInj}) gives
        \[
        ||\sigma(x) - \sigma(y)|| \geqslant (\min\{C,C_+,C_-\}) \cdot
        3^{-\beta(1-\zeta)} \cdot |\gamma|^{\beta(1-\zeta)} \ .
        \]
        Hence, there exists $C'_\zeta >0$ with $||\sigma(x) - \sigma(y)||
        \geqslant C'_\zeta \cdot |x^{-1}y|^{\beta(1-\zeta)}$.
    \end{itemize}
    Consequently, setting $m_\zeta = \min \left\{ \alpha\zeta , \beta(1-\zeta) \right\}$, it comes
    \begin{equation*}\label{eqWr2}\tag{$\ast\ast_\zeta$}
    || \sigma(x) - \sigma(y) || \geqslant \min \left\{ \frac C{3^{\alpha\zeta}}, C'_\zeta \right\} \cdot |\gamma|^{m_\zeta}  \ .
    \end{equation*}
    The largest value for $m_\zeta$ is obtained for $\alpha\zeta = \beta(1-\zeta)$, that is
    $\zeta = \frac \beta{\alpha + \beta}$. It gives $m_\zeta = \frac {\alpha\beta}{\alpha+\beta}$.
    This is our second estimate for case (c).
\end{description}

As one has $\beta > \alpha\beta/(\alpha+\beta)$, combination of cases (a)-(c) gives
\begin{align*}
|| \sigma(x) - \sigma(y) || &\geqslant C'' \cdot |x^{-1} y|^{\frac{\alpha\beta}{\alpha + \beta}}
& \forall x,y \in \Gamma \\
|| \sigma(x) - \sigma(y) || &\geqslant C'' \cdot |x^{-1} y|^{\min\{\beta, \alpha - \frac 12\}}
& \forall x,y \in \Gamma
\end{align*}
for some $C'' >0$.
Hence, we get $R_\sigma \geqslant \alpha\beta/(\alpha+\beta)$ and $R_\sigma \geqslant
\min\{\beta, \alpha - \frac 12\}$ for all $\alpha,\beta$ satisfying $0<\alpha<R_f$
and $0<\beta<R_{f_\pm}$. This implies immediately $R_\sigma \geqslant R_f\cdot R_{f_\pm}/(R_f+R_{f_\pm})$
and $R_\sigma \geqslant \min\{R_{f_\pm}, R_f - \frac 12\}$.

To conclude the proof of Proposition \ref{HilbWr}, we pass now to the last statement. We thus suppose that $f$ is $H$-equivariant and
$f_\pm$ are $\Gamma$-equivariant (with respect to some actions by affine isometries). To establish the
$\Gamma$-equivariance of $\sigma$, we only have to define a $\Gamma$-action (by affine isometries) on
$\oplus_{i\in\Z} \mathcal H$ and check the $\Gamma$-equivariance with respect to it.

The $\Gamma$-action on $\bigoplus_{i\in\Z}\mathcal{H} = \ell^2(\Z,\mathcal{H},0)$ is defined by equation
(\ref{ActWr}). To check the equivariance, we set $\gamma = (\lambda,n)$
and $g=(\mu,p)$ with ${\lambda,\mu \in H^{(\Z)}}$ and $n,p \in\Z$. We have $(\gamma \cdot \sigma(g))_i =
\lambda_i \cdot f(\mu_{i-n})$ and $\sigma(\gamma g)_i = f(\lambda_i \mu_{i-n})$ and we get $(\gamma \cdot
\sigma(g))_i = \sigma(\gamma g)_i$ for all $i$ by $H$-equivariance of $f$.
\end{proof}

\medskip

Theorem \ref{HilbC} will be obtained by applying Proposition \ref{HilbWr} with good embeddings of the trees
$T_\pm$. We explain now how to embed a tree in a Hilbert space with high values of the constant ``$R_f$''.
First, the following result can be obtained by a straightforward adaptation of \cite[Proposition 4.2]{GuKa}.
\begin{Prop}\label{CompTree}
Let $T= (V,E)$ be a tree. Then $R(V) = 1$.
\end{Prop}
More precisely, if we denote by $E_G$ the set of geometric (or unoriented) edges of $T$ and if we fix a
base vertex $v_0$, then for any $\varepsilon \in \ ]0,1/2[ \ $ we may consider the map
\[
f_\varepsilon \ : \ V \longrightarrow \ell^2(E_G) \ ; \ x \longmapsto \sum_{k=1}^{d(v_0,x)} k^\varepsilon \delta_{e_k(x)} \ ,
\]
where the $e_k(x)$'s are the consecutive edges on the unique geodesic from $x$ to $v_0$ and $\delta_e$ is the Dirac mass
at $e$. It is a Lipschitz map with $R_{f_\varepsilon} \geqslant 1/2 + \varepsilon$. We refer to the proof of
\cite[Proposition 4.2]{GuKa} for this fact.

To prove the ``equivariant'' part of Theorem \ref{HilbC}, we need some explicit equivariant embeddings into
Hilbert spaces. Let $T=(V,E)$ be a tree. We recall from Section 7.4.1
in \cite{CCJJV} how to embed equivariantly $T$ in a
Hilbert space.
%Let $\mathbf{E}$ be the set of oriented edges of $T$.
We recall that we denote by $e\mapsto\overline{e}$
the ``orientation-reversing'' involution on $E$, and we endow
$\ell^2(E)$ with the scalar product:
$$\langle\xi|\eta\rangle = \frac{1}{2} \sum_{e\in E} \xi(e)\overline{\eta(e)}.$$
Define a map $c:V\times V\rightarrow\ell^2(E):(x,y)\mapsto c(x,y)$
with
$$c(x,y)=\sum_{e\in(x\rightarrow y)}\delta_e -\delta_{\overline{e}}$$
where $\delta_e$ is the Dirac mass at $e$ and the summation is taken over coherently oriented
edges in the oriented geodesic from $x$ to $y$. The map $c$ satisfies, for every $x,y,z \in V$:
\begin{equation}\label{eq1EmTrH}
c(x,y)+c(y,z)=c(x,z);
\end{equation}
\begin{equation}\label{eq2EmTrH}
\|c(x,y)\|^2 = d(x,y).
\end{equation}
Moreover if a group $G$ acts on $T$, then for every $g\in G$:
\begin{equation}\label{eq3EmTrH}
c(gx,gy)=\pi(g)c(x,y)
\end{equation}
where $\pi$ is the permutation representation of $G$ on $\ell^2(E)$.

Fix now a base-vertex $v_0\in V$. Define a map
\[
\iota_{T,v_0}:V\rightarrow \ell^2(E):v\mapsto c(v_0,v)
\]
and, for $g\in G$, an affine isometry $\alpha_{v_0}(g)$ of $\ell^2(E)$:
\[
\alpha_{v_0}(g)\xi =\pi(g)\xi + c(v_0,gv_0).
\]
Using equations (\ref{eq1EmTrH}) -- (\ref{eq3EmTrH}) above, the following lemma is immediate.

\begin{Lem}\label{equiv}
\begin{enumerate}
\item For all $g,h\in G:\,\alpha_{v_0}(gh)=\alpha_{v_0}(g)\alpha_{v_0}(h)$, so that $\alpha_{v_0}$ defines an affine isometric action
of $G$ on $\ell^2(E)$;
\item the map $\iota_{T,v_0}$ is $G$-equivariant with respect to the action $\alpha_{v_0}$ on $\ell^2(E)$;
\item one has $|| \iota_{T,v_0}(x) - \iota_{T,v_0}(y) || = \sqrt{d(x,y)}$ for all $x,y\in V$, so that $R_{\iota_{T,v_0}} =
1/2$.
\end{enumerate}
\end{Lem}
It is an immediate consequence that $R_G(V) \geqslant 1/2$.

\begin{sproof}{Theorem \ref{HilbC}}
The inequalities $R(H) \geqslant R(\Gamma)$ and $R_H(H) \geqslant R_\Gamma(\Gamma)$ are trivial.

One has $R(V(T_\pm)) = 1$ by Proposition \ref{CompTree}, so that Proposition \ref{HilbWr} gives
$R(\Gamma) \geqslant R(V(T_\pm))\cdot R(H)/(R(V(T_\pm)) + R(H)) = R(H)/(R(H)+1)$.

Finally, one has $R_\Gamma(V(T_\pm)) \geqslant 1/2$ by Lemma \ref{equiv}, so that we obtain
\begin{align*}
R_\Gamma(\Gamma) &\geqslant \frac{R_\Gamma(V({T_\pm}))\cdot R_H(H)}{R_\Gamma(V(T_\pm)) + R_H(H)} \geqslant
\frac{R_H(H)}{2R_H(H)+1} \\
R_\Gamma(\Gamma) & \geqslant \min \left\{ R_\Gamma(V(T_\pm)) , R_H(H)- \frac 12 \right\} = R_H(H) - \frac 12
\end{align*}
by Proposition \ref{HilbWr}.
\end{sproof}

%---------------------------------------------------------------------------------------------------
\section{Hilbert space compression: examples}\label{Ex}
We begin this section with known results about the compression of groups of the form $H\wr\Z$.
Let us first state a generalization of \cite[Theorem 3.9]{AGS} which gives upper bounds for many
of them.

\begin{Prop}\label{upC}
Let $G$ be a finitely generated group with growth function satisfying $\kappa(n) \succcurlyeq n^k$ for some $k>0$
and let $H$ be a group. We assume the generating set of $H$ chosen such that the word metric is unbounded.
Then, the Hilbert space compression of $\Gamma = H\wr G$ satisfies
\[
R(\Gamma,\Sigma) \leqslant \frac{1+ k/2}{1+k} \ ,
\]
where $\Sigma$ is the union of the generating sets of $G$ and $H$. In particular, with $G = \Z$, we get
$R(H\wr\Z)\leqslant 3/4$.
\end{Prop}
The proof is a straightforward adaptation of \cite[Theorem 3.9]{AGS}.
\begin{Rem}
If $H$ is finitely generated, the hypothesis ``the word metric is unbounded'' means exactly that $H$ is
infinite.
\end{Rem}
Lower bounds, were found by Tessera \cite[Corollary 14]{Tess}. In particular:
\begin{Prop}\label{lowC}
Let $H$ be a finitely generated group. If $H$ has polynomial growth, one has $R(H\wr\Z) \geqslant 2/3$.
\end{Prop}
Together, Propositions \ref{upC} and \ref{lowC} give immediately:
\begin{Cor}\label{Cor1}
If $H$ is an infinite group with polynomial growth, then one has $R(H\wr\Z)\in [2/3,3/4]$.
\end{Cor}

%Arzhantseva, Guba and Sapir had originally proved that $R(\Z\wr\Z)\geqslant \frac{1}{2}$
%\cite[Theorem 1.8]{AGS}, using the notion of diagram group. We observe that this original
%lower bound can be obtained without using diagram groups. It suffices to apply Theorem
%\ref{HilbC} since the Hilbert space compression of $\Z$ is $1$.

%The other famous example of group with Hilbert space compression different from $0$ and $1$ is the following:
%\begin{Prop}\label{Cor0}\cite[Theorems 1.3 and 1.13]{AGS}
%Thompson's group $F$ satisfies $R(F) = 1/2 = R_F(F)$.
%\end{Prop}

%We now provide new examples of such groups, whose
%construction can obviously be iterated.

%\begin{Cor}\label{Cor2}
%\begin{enumerate}
%  \item The group $F\wr\Z$ (where $F$ is Thompson's group) satisfies $R(F\wr\Z) \in [1/3,1/2]$ and
%  $R_{F\wr\Z}(F\wr\Z) \in [1/4,1/2]$.
% \item One has ${R((\Z\wr\Z)\wr\Z),R_{(\Z\wr\Z)\wr\Z} ((\Z\wr\Z)\wr\Z) \in [2/5,3/4]}$;
%\end{enumerate}
%\end{Cor}
%\begin{proof}
%Part 1 follows from Theorem \ref{HilbC} and Proposition \ref{Cor0} and part 2 from Theorem \ref{HilbC},
%Corollary \ref{HilbCmoy} and Corollary \ref{Cor1}.
%(We notice here that Theorem \ref{HilbC} would only give $R_{(\Z\wr\Z)\wr\Z} ((\Z\wr\Z)\wr\Z) \in [2/7,3/4]$.)
%\end{proof}

%\medskip

%The exact value of $R_{F\wr\Z}(F\wr\Z)$ could be of interest, since the inequality $R_{F\wr\Z}(F\wr\Z) < 1/3$
%would imply that $F$ is non amenable.

In a similar spirit, Proposition \ref{upC} and our Theorem \ref{HilbC} imply immediately:

\begin{Cor}\label{Examp} Let $H$ be an infinite, finitely generated group.
\begin{enumerate}
\item [a)] If $R(H)=1$, then $R(H\wr\Z)\in [1/2,3/4]$.
\item [b)] If $R(H)=R_H(H)=1/2$, then $R(H\wr\Z)\in [1/3,1/2]$ and $R_{H\wr\Z}(H\wr\Z)\in [1/4,1/2]$ (in particular, if
$R_{H\wr\Z}(H\wr\Z)<1/3$, then $H$ is non-amenable).

\hfill$\square$
\end{enumerate}
\end{Cor}

The interest of part (a) in Corollary \ref{Examp} stems from the fact that numerous groups satisfy $R(H)=1$: among amenable groups,
we mention polycyclic groups and lamplighter groups $F\wr\Z$ with $F$ finite \cite[Theorem 1]{Tess}; among (usually) non-amenable
groups, we cite hyperbolic groups \cite[Theorem 4.2]{BrodSon}, groups acting properly co-compactly on finite-dimensional $CAT(0)$
cube complexes \cite{CampNib}, co-compact lattices in connected Lie groups, irreducible lattices in higher rank semi-simple Lie groups
\cite[Theorem 2]{Tess}.

Our excuse for isolating (b) in Corollary \ref{Examp} is a remarkable result by Arjantseva, Guba and Sapir \cite[Theorem 1.8]{AGS}:
for Thompson's group $F$, one has $R(F)=R_F(F)=1/2$.

\paragraph{Acknowledgements.} The authors thank Yves de Cornulier for his relevant remarks on a previous version.

\medskip

%-------------------------------------------------------------------------------------------------

\bibliographystyle{alpha}

%\bibliography{MaBiblio}

\newcommand{\etalchar}[1]{$^{#1}$}
\def\cprime{$'$}

%-----------------------------------------------------------------------------

\medskip

Authors address: \vspace{2mm} \\
Institut de Math\'ematiques \\
Universit\'e de Neuch\^atel \\
Rue Emile Argand 11 \\
Case postale 158 \\
CH-2009 Neuch\^atel \\
SWITZERLAND

\vspace{2mm}

yves.stalder@unine.ch; alain.valette@unine.ch

\end{document}